\documentclass{article}
\usepackage{amsmath,amsfonts,latexsym,amssymb,euscript,amsthm}

\newcommand{\Om}{\Omega}

\renewcommand{\d}{\partial}

\newcommand{\ep}{\varepsilon}
\newcommand{\intl}{\int\limits}
\renewcommand{\div}{\mathop{\mathrm{div }}}
\newcommand{\rot}{\mathop{\mathrm{rot }}}

\newtheorem{theorem}{Theorem}[section]

\newtheorem{lemma}{Lemma}[section]
\newtheorem{corollary}{Corollary}[section]
\renewcommand{\l}{\left}
\renewcommand{\r}{\right}

\newcommand{\Q}{Q^+}
\newcommand{\B}{B^+}
\newcommand{\al}{\alpha}

\newcommand{\gr}{\nabla}

\renewcommand{\div}{\mathop{\mathrm{div }}}
\newcommand{\Ga}{\Gamma}

\renewcommand{\L}{\EuScript{L}}

\newcommand{\tht}{\theta}
\newcommand{\cd}{\partial}

\title{On the local smoothness of weak solutions to the MHD system near the boundary}
\author{V.~Vyalov
\footnote{This research is supported by the Chebyshev Laboratory
(Department of Mathematics and Mechanics, St.-Petersburg State University)
under RF government grant 11.G34.31.0026}
}

\begin{document}

\maketitle

\section{Introduction}

Assume $\Om\subset \mathbb R^3$ is a $C^2$-- smooth bounded domain
and $Q_T=\Om\times (0,T)$. In this paper we investigate the boundary
regularity of solutions to the principal system of
magnetohydrodynamics (the MHD equations):
\begin{equation}
\left.
\begin{array}c
\partial_t v  + (v\cdot \nabla )v  - \Delta v +  \nabla
p = \rot H \times H
\\
 \div v =0
\end{array} \right\}
\quad\mbox{in } Q_T, \label{MHD_NSE}
\end{equation}
\begin{equation}
\left.
\begin{array}c
\partial_t H    +\rot \rot H = \rot (v\times H)
\\ \div H =0
\end{array} \right\} \quad \mbox{in  }Q_T.
\label{MHD_Magnetic}
\end{equation}
Here  unknowns
are  the velocity field $v:Q_T\to \mathbb R^3$, pressure
$p:Q_T\to \mathbb R$, and the magnetic field $H:Q_T\to \mathbb
R^3$. We impose on $v$ and $H$ the boundary conditions:
\begin{equation}
v|_{\cd\Om\times(0,T)}=0, \quad H_{\nu}|_{\cd\Om\times(0,T)}=0, \quad (\rot H)_{\tau}|_{\cd\Om\times(0,T)}=0,
\label{BC}
\end{equation}
Here by $\nu$ we denote the outer normal to $\cd\Om$ and $H_{\nu}= H\cdot \nu$, $(\rot H)_\tau=\rot H- \nu(\rot H\cdot \nu)$.
We introduce the following definition:

\medskip\noindent{\bf Definition:} Assume $\Ga\subset \cd\Om$.
 The functions $(v,H, p)$ are called a {\it boundary suitable
 weak solution } to the system  (\ref{MHD_NSE}), (\ref{MHD_Magnetic})
near $\Ga_T\equiv \Ga\times (0,T)$ if
\begin{itemize}
\item[1)] \ $v\in L_{2,\infty}(Q_T)\cap W^{1,0}_2(Q_T)\cap W^{2,1}_{\frac 98, \frac 32}(Q_T)$,

$H \in L_{2,\infty}(Q_T)\cap W^{1,0}_2(Q_T)$,
\item[2)] \ $p\in L_{\frac 32}(Q_T)\cap W^{1,0}_{\frac 98, \frac 32}(Q_T)$,
\item[3)] \ $\div v=0$, \ $\div H=0$ \ a.e. in \ $Q_T$,
\item[4)] \ $v|_{\cd \Om}=0$, \ $H_\nu|_{\cd\Om}=0$ \ in the sense of traces,
\item[5)] \ for any $w\in L_2(\Om)$ the functions
$$
t\mapsto \int\limits_{\Om} v(x,t)\cdot w(x)~dx \qquad\mbox{and}\qquad t\mapsto \int\limits_{\Om} H(x,t)\cdot w(x)~dx
$$
are continuous,

\item[6)] \ $(v,H)$ satisfy the following integral identities: for any $t\in [0,T]$
$$
\begin{array}c
\int\limits_{\Om} ~v(x,t)\cdot \eta(x,t)~dx -\int\limits_{\Om} ~v_0(x)\cdot \eta(x,0)~dx \ + \\ + \
\int\limits_0^t \int\limits_\Om ~ \Big( - v\cdot \cd_t \eta + (\nabla v  - v\otimes v + H \otimes H) : \nabla \eta - (p+\frac 12 |H|^2)\div \eta \Big)~dxdt \ = \ 0,
\end{array}
$$
for all \ $\eta \in W^{1,1}_{\frac 52}(Q_t)$ \ such that   \ $\eta|_{\cd\Om\times (0,t)}=0$,
$$
\begin{array}c
 \int\limits_{\Om}~ H(x,t)\cdot \psi (x,t)~dx  \ -  \  \int\limits_{\Om}~ H_0(x)\cdot \psi (x,0)~dx \ + \\ + \
\int\limits_0^t \int\limits_\Om ~ \Big( - H\cdot \cd_t \psi  + \rot H\cdot \rot \psi - (v\times H) \cdot \rot \psi \Big)~dxdt  \ = 0,
\end{array}
$$
for all \ $\psi \in W^{1,1}_{\frac 52}(Q_t)$ \ such that \  $\psi_\nu |_{\cd\Om\times (0,t)}=0$.

\item[7)] For every $z_0=(x_0,t_0)\in \Ga_T$  such that  $ \Om_R(x_0)\times (t_0-R^2, t_0)\subset Q_T$
and for any $\zeta\in C_0^\infty(B_R(x_0)\times(t_0-R^2,t_0])$ such that $\left.\frac{\cd\zeta}{\cd \nu}\right|_{\cd\Om}=0$
the following ``local energy inequality near $\Ga_T$'' holds:
\begin{equation}
\begin{array}c
\sup\limits_{t\in (t_0-R^2,t_0)}\int\limits_{\Om_R(x_0)}  \zeta \Big( |v|^2+|H|^2\Big)~dx \ + \\
 + \ 2 \int\limits_{t_0-R^2}^{t_0}\int\limits_{\Om_R(x_0)} \zeta \Big( |\nabla v|^2+|\rot H|^2\Big)~dxdt \ \le \\
\le\  \int\limits_{t_0-R^2}^{t_0}\int\limits_{\Om_R(x_0)}  \Big(| v|^2 +|H|^2\Big) (\cd_t \zeta + \Delta \zeta )~dxdt \ + \\ + \
\int\limits_{t_0-R^2}^{t_0}\int\limits_{\Om_R(x_0)} \Big ( | v|^2+ 2\bar p \Big ) v\cdot\nabla \zeta~dxdt \ + \\
- \ 2 \int\limits_{t_0-R^2}^{t_0}\int\limits_{\Om_R(x_0)}   (H\otimes H):  \nabla^2 \zeta  ~dxdt \ + \\ + \
2 \int\limits_{t_0-R^2}^{t_0}\int\limits_{\Om_R(x_0)}  (v\times  H)(\nabla \zeta\times H)~dxdt
\end{array}
\label{LEI}
\end{equation}
We remark also that the following identity holds
$$
\begin{array}c
 \int\limits_{t_0-R^2}^{t_0}\int\limits_{\Om_R(x_0)}  (v\times  H)(\nabla \zeta\times H)~dxdt  \ =  \\ = \
\int\limits_{t_0-R^2}^{t_0}\int\limits_{\Om_R(x_0)} ~ (v \cdot \nabla \zeta) |H|^2  ~dxdt  \ -
\ \int\limits_{t_0-R^2}^{t_0}\int\limits_{\Om_R(x_0)} ~ (v \cdot H) (H\cdot \nabla \zeta)  ~dxdt
\end{array}
$$
\end{itemize}

\noindent
Here  $L_{s,l}(Q_T)$ is the anisotropic Lebesgue space equipped with the norm
$$
\|f\|_{L_{s,l}(Q_T)}:=
\Big(\int_0^T\Big(\int_\Om |f(x,t)|^s~dx\Big)^{l/s}dt\Big)^{1/l} ,
$$
and we use the following notation for the functional spaces:
$$
\gathered
W^{1,0}_{s,l}(Q_T)\equiv L_l(0,T; W^1_s(\Om))= \{ \ u\in L_{s,l}(Q_T): ~\nabla u \in L_{s,l}(Q_T) \ \},
\\
W^{2,1}_{s,l}(Q_T) = \{ \ u\in W^{1,0}_{s,l}(Q_T): ~\nabla^2 u, \ \cd_t u \in L_{s,l}(Q_T) \ \},
\\
\overset{\circ}{W}{^1_s}(\Om)=\{ \ u\in W^1_s(\Om):~ u|_{\cd\Om}=0 \ \},
\endgathered
$$
and the following notation for the norms:
$$
\gathered
\| u \|_{W^{1,0}_{s,l}(Q_T)}= \| u \|_{L_{s,l}(Q_T)}+ \|\nabla u\|_{L_{s,l}(Q_T)},  \\
\| u \|_{W^{2,1}_{s,l}(Q_T)}= \| u \|_{W^{1,0}_{s,l}(Q_T)}+ \| \nabla^2 u \|_{L_{s,l}(Q_T)}+\|\cd_t u\|_{L_{s,l}(Q_T)},  \\
\endgathered
$$

Denote $B(x_0,r)$ the open ball in $\mathbb R^3$ of radius $r$ centered at $x_0$ and denote by $B^+(x_0,r)$ the half--ball $\{ x\in B(x_0,r)~|~ x_3>0\}$.
For $z_0=(x_0,t_0)$ denote $Q(z_0,r)=B(x_0,r)\times (t_0-r^2, t_0)$, $Q^+(z_0,r)=B^+(x_0,r)\times (t_0-r^2, t_0)$.
In this paper we shall use the abbreviations: $B(r)= B(0,r)$, $B^+(r)= B^+(0,r)$ etc, $B=B(1)$, $B^+=B^+(1)$ etc.

\section{Main Results}

Our work deals with the sufficient conditions  of local regularity of suitable weak solutions to the MHD system near the plane part of the boundary.
In \cite{VyaShi} the following results were obtained.
\begin{theorem}
\label{Main_ep-reg}
There exists an absolute constant $\ep_*>0$ with the following property. Assume $(v,H,p)$ is a boundary suitable weak solution in $Q_T$
and assume $z_0=(x_0,t_0)\in \d\Om\times (0,T)$ is such that $x_0$ belongs to the plane part of $\d\Om$.
If there exists $r_0>0$ such that $Q^+(z_0,r_0)\subset Q_T$ and
$$
\frac 1{r_0^2} \int\limits_{Q^+(z_0,r_0)} \Big(~ |v|^3+|H|^3+|p|^{\frac 32}~\Big) dxdt \  < \ \ep_*,
$$
then the functions $v$ and $H$ are H\" older continuous on $\bar Q^+(z_0, \frac {r_0}2)$.
\end{theorem}

\begin{theorem}\label{CKN_theorem}
For any  $K>0$ there exists $\ep_0(K)>0$ with the following property.
Assume $(v,H,p)$ is a boundary suitable weak solution in $Q_T$ and assume $z_0=(x_0,t_0)\in \cd\Om\times (0,T)$ is such that $x_0$ belongs to the plane part of $\cd\Om$.
If
\begin{equation}
\limsup\limits_{r\to 0}\Big(~\frac 1r \int\limits_{Q(z_0,r)} |\nabla H|^2~dxdt ~\Big)^{1/2}\ < \ K
\label{ep-regularity-1}
\end{equation}
and
\begin{equation}
\limsup\limits_{r\to 0} \Big(~\frac 1r \int\limits_{Q(z_0,r)} |\nabla v|^2~dxdt~\Big)^{1/2} \ < \ \ep_0,
\label{ep-regularity-2}
\end{equation}
then there exists $\rho_*>0$ such that the functions $v$ and $H$ are H\" older continuous on the closure of $ Q^+(z_0,\rho_*)$.
\end{theorem}


Let us explain  the main differences between these theorems. The statement of the theorem \ref{Main_ep-reg} contains smallness
conditions on the three functionals, but these conditions have to hold only for one value of cylinder radius.
In theorem \ref{CKN_theorem} we have conditions for all sufficiently small values of radius, but smallness condition \eqref{ep-regularity-2}
is imposed only on the velocity $v$.

To describe more conditions of local regularity we will need the following notations
\begin{equation}
\begin{array}c
E(r) \ = \
\Big( ~\frac 1{r} \int\limits_{Q^+(r)} |\nabla v|^2~dxdt ~ \Big)^{1/2},
\\
 E_*(r) \ = \  \Big( \frac 1{r} \int\limits_{Q^+(r)} |\nabla H|^2~dxdt \Big)^{1/2},
\\
A( r) \equiv \Big( \frac 1{r}
\sup\limits_{t\in (-r^2, 0)} \int\limits_{B^+(r)} |v|^{2}~dy
\Big)^{1/2}, \\
A_*( r) \equiv \Big( \frac 1{r}
\sup\limits_{t\in (-r^2, 0)} \int\limits_{B^+(r)} |H |^{2}~dy
\Big)^{1/2},
\\
C_q( r) \equiv \Big( \frac 1{r^{5-q}} \int\limits_{Q^+(r)} |v|^q~dydt
\Big)^{1/q},
\\
F_q( r) \ = \ \Big( \frac 1{r^{5-q}} \int\limits_{Q^+(r)} |H|^q~dxdt
\Big)^{1/q}
\\ D(r) \equiv \Big( \frac 1{r^2} \int\limits_{Q^+(r)}
|p - [p]_{B^+(r)}|^{3/2}~dydt \Big)^{2/3},
\\
D_s(r) = R^{\frac 53 -\frac 3s} \Big( \int\limits_{-r^2}^{0} \Big(
\int\limits_{B^+(r)}  |\nabla p|^{s}~dy \Big)^{\frac 1s \cdot
\frac 32}~dt \Big)^{2/3},
\end{array}
\label{EF}
\end{equation}
$$
C(r) = C_3(r), \qquad F(r)= F_3(r), \qquad D_*(r)= D_{\frac{36}{35}}(r).
$$

Note that the equations \eqref{MHD_NSE},  \eqref{MHD_Magnetic}, as well as the functionals \eqref{EF}
and the statements of the previous theorems are invariant under the scaling transformations
\begin{equation}
\left.
\begin{aligned}
& v_{\rho}(y,s) = \rho v(\rho y + x_0, \rho^2 s + t_0),\\
& H_{\rho}(y,s) = \rho H(\rho y + x_0, \rho^2 s + t_0),\\
& p_{\rho}(y,s) = \rho^2 p(\rho y + x_0, \rho^2 s + t_0).
\label{scaling}
\end{aligned}
\right\}
\end{equation}

We use the approach which was originally  developed in \cite{Ser1} for the Navier-Stokes equations (and later  it was used also in \cite{Mih}). According to this approach the regularity of solutions follows if  one of the functionals \eqref{EF} is bounded  uniformly with respect to $r$ and additionally one of these functionals is small only for a single  sufficiently small value of the radius.
Our goal is to obtain the same result for the solutions to the MHD system.

The main result of our work is the following theorem, that is a kind of ``interpolation'' of theorems \ref{Main_ep-reg} and \ref{CKN_theorem}.
\begin{theorem}\label{Seregin_theorem}
For arbitrary  $K>0$ there is a constant $\ep_1(K)>0$ with the following property.
Assumw $(v,H,p)$ is a suitable weak solution to the MHD system in $Q_T$ and $z_0=(x_0,t_0)\in \d\Om\times (0,T)$ where $x_0$
belongs to the plane part of $\d\Om$.
If
\begin{equation}
\limsup\limits_{r\to 0}\Big(~\frac{1}{r^2} \int\limits_{\Q(z_0,r)}|v|^3~dxdt ~\Big)^{1/3} +
 \Big(~\frac{1}{r^3} \int\limits_{\Q(z_0,r)}|H|^2~dxdt ~\Big)^{1/2}\ < \ K
\end{equation}
and one of the following conditions holds
\begin{equation}
\gathered
\liminf\limits_{r\to 0}\Big(~\frac 1r \int\limits_{\Q(z_0,r)} |\nabla v|^2~dxdt~\Big)^{1/2} < \ \ep_1, \\
\liminf\limits_{r\to 0}\Big(~\frac 1r \sup\limits_{-r^2<t<0} \int\limits_{\B(x_0,r)} |v|^2~dxdt~\Big)^{1/2} \ < \ \ep_1, \\
\liminf\limits_{r\to 0}\Big(~\frac{1}{r^2} \int\limits_{\Q(z_0,r)} |v|^3~dxdt~\Big)^{1/3} < \ \ep_1,
\endgathered
\end{equation}
then there exists $\rho_*>0$ such that the functions $v$ and $H$ are H\" older continuous on the closure of $ Q^+(z_0,\rho_*)$.
\end{theorem}

Note that it is possible to prove a lot of analogues of theorem \ref{Seregin_theorem}. Generally the proof consist of two steps.
The first step  is the proof of boundedness of energy functionals \eqref{EF}. Usually to do this it is sufficient to have boundedness condition
for one functional depending on $v$ and for another one depending on $H$.  In our work this step is carried out in  Section \ref{Bound_section}. Also
we will use some estimates for the magnetic field $H$, that can be obtained if we consider equation \eqref{MHD_Magnetic} as the heat equation with lower order terms depending on $v$. This inequalities are proved in Section \ref{Hdecomposition_section}.

The second step is the proof of regularity condition when  all of functionals \eqref{EF} are bounded and one of functionals on $v$ is small for a single  sufficiently small value of $r$.
This result can be found in Section \ref{reg_section}.


\section{Estimates of Solutions to the Heat Equation}
\label{Hdecomposition_section}

In this section we study solutions of the heat equations  with the lower order terms:
$$
\gathered
\d_t H \ - \ \Delta H \ =  \ \div (v\otimes   H - H \otimes v) \qquad\mbox{in}\quad Q^+.
\\
v|_{x_3=0}=0,
\\
H_3|_{x_3=0}=0, \qquad  H_{\al, 3} |_{x_3=0}=0, \qquad \al=1,2.
\endgathered
$$
Namely, we assume the functions $(v, H)$ possess the following properties:
\begin{equation}
\gathered
v, \ H\in W^{1,0}_2(Q^+), \\ v|_{x_3=0}=0, \quad H_3|_{x_3=0}=0 \quad \mbox{in the sense of traces}, \\
\endgathered
\label{Class}
\end{equation}
for any $\eta \in C^\infty_0(Q; \mathbb R^3)$ such that  $\eta_3|_{x_3=0}=0$ the following integral identity holds
\begin{equation}
\gathered
\int\limits_{Q^+} \Big( - H\cdot \d_t \eta + \nabla H:\nabla \eta \Big)~dxdt \ =  \ - \int\limits_{Q^+} G:\nabla \eta~dxdt,  \\
\endgathered
\label{Integral_Identity}
\end{equation}
here $G\ = \ v\otimes H - H\otimes v$, and
\begin{equation}
\div v=0, \qquad\div H=0 \qquad\mbox{a.e. in}\quad Q^+.
\label{Divergent-free}
\end{equation}

\begin{lemma}
Assume that conditions (\ref{Class}) --- (\ref{Divergent-free}) hold. Then for any
$0< r \leq 1$ and $0 < \tht \leq 1$ the following estimate holds
\begin{equation}
F_2(\tht r) \leq c \tht^{\al} F_2(r) + c \tht^{-\frac32} C(r) A_*(r).
\label{decHC}
\end{equation}
\end{lemma}

\bigskip
\noindent
{\bf Proof.}
Denote by $v^*$ and $H^*$ the extensions of functions $v$ and $H$ from $Q^+$ onto $Q$.
Fix  arbitrary $r\in (0, 1)$ and let $\zeta\in C^\infty(\bar Q)$ be a cut off function such that $\zeta \equiv 1$ on $Q(r)$ and
$\operatorname{supp} \zeta\subset B\times (-1,0]$. Denote $ \Pi=\mathbb R^3\times (-1,0)$ and denote by ${\hat G}$ the function which coincides with
$G^*$ on $Q(\frac r2)$ and additionally possesses the following properties: ${\hat G}\in W^{1,0}_1(\Pi)\cap L_{\frac {18}{11}, \frac 65}(\Pi)$,
${\hat G}$ is compactly supported in $\Pi$,  and
\begin{equation}
\| {\hat G}\|_{L_{\frac65, 2}(\Pi)} \ \le \ c \| G^{*}\|_{L_{\frac65,2}(Q(\frac r2))} \ \le \ c \| G\|_{L_{\frac65,2}(Q^+(\frac r2))}
\label{G**}
\end{equation}

We decompose $H^*$ as
$$
H^* \ =  \hat H \ + \tilde H,
$$
where $\hat H$ is a solution of the  Cauchy problem for the heat equation
\begin{equation}
\left\{ \ \begin{array}l
\d _t \hat H -\Delta \hat H  \ = \ \div {\hat G} \qquad \mbox{in}\quad \Pi, \\
\hat H |_{t=-1} =0,
\end{array}\right.
\label{Cauchy problem}
\end{equation}
defined by the formula
$\hat H = \Ga * \div  {\hat G}= -\nabla \Ga *  {\hat G}$, where $\Ga$ is the fundamental solution of the heat operator.
The function $\tilde H$ satisfies the homogeneous heat equation
\begin{equation}
\begin{array}c
\d_t \tilde H - \Delta \tilde H =  0 \qquad \mbox{in}\quad Q( \frac r2).
\end{array}
\label{Heat-2}
\end{equation}

Take arbitrary $\theta \in (0,\frac 12)$.
We estimate $\| H\|_{L_2(Q^+(\theta r))}$ in the following way
\begin{equation}
\gathered
\| H\|_{L_2(Q^+(\theta r))} \ \le \ \| H^* \|_{L_2(Q(\theta r))}  \ \le
\ \| \hat H\|_{L_2(Q(\theta r))} \ + \ \| \tilde H\|_{L_2(Q(\theta r))},
\endgathered
\label{F-2}
\end{equation}
For  $\| \hat H\|_{L_2(Q(\theta r))}$ we have
\begin{equation}
\| \hat H\|_{L_2(Q(\theta r))} \ \le \ c~\| \hat H\|_{L_2(Q(\frac  r2))}.
\label{F-3}
\end{equation}
As $\tilde H$ satisfies (\ref{Heat-2}) by local estimate of the maximum of $\tilde H$ via its $L_2$--norm  we obtain
\begin{equation}
\gathered
\| \tilde H\|_{L_2(Q(\theta r))} \ \le \ c ~\theta^{\frac {5}{2}} ~\| \tilde H\|_{L_2(Q(\frac r2))} \ \le
\\ \le
\ c ~\theta^{\frac {5}{2} } ~( \|  H^*\|_{L_2(Q(r))} + \| \hat H\|_{L_2(Q(\frac r2))})
\endgathered
\label{F-4}
\end{equation}

So, we need to estimate $\| \hat H\|_{L_2(Q(\frac r2))}$. As singular integrals are bounded
on the anisotropic Lesbegue space $L_{s,l}$ (see, for example, \cite{Solonnikov_Uspekhi}) for the convolution  $\hat h = \Ga* {\hat G} $ we obtain the estimate
$$
\| \hat h \|_{W^{2,1}_{\frac65, 2}(Q(r))} \ \le  \ c \|  {\hat G} \|_{L_{\frac65,2}(\Pi)}.
$$
On the other hand, from  the 3D-- parabolic imbedding theorem (see \cite{Besov})
$$
W^{2,1}_{s,l}(Q)\hookrightarrow W^{1,0}_{p,q}(Q), \quad\mbox{as}\quad 1-\left(\frac 3s+\frac 2l - \frac 3p -\frac 2q\right) \ge 0,
$$
for $p=q=2$ and $s=\frac65$, $l=2$ and for $\hat H=- \nabla \hat h$ we obtain
$$
\| \hat H \|_{L_2(Q(r))} \ \le \ c ~\|  {\hat G} \|_{L_{\frac65, 2}(\Pi)}.
$$
(Note that the constant $c$ in this inequality does not depend on $r$). Taking into account (\ref{G**}) we arrive at
\begin{equation}
\| \hat H \|_{ L_2(Q(r))} \ \le \ c ~ \| G \|_{L_{\frac{6}{5}, 2}(Q^+(\frac r2))}.
\label{imbedding}
\end{equation}

From the definition of $G$ we obtain
\begin{equation}
\gathered
\| G \|_{L_{\frac65, 2}(Q^+(\frac r2))} \ \le \ c ~
\left(\ \int\limits_{-r^2/4}^{0} \| v\otimes H\|_{L_{\frac65}(B^+(r/2))}^{2}~dt\ \right)^{\frac12} \leq \\
\leq \l( \intl_{-r^2}^0 \| v \|_{3,\B(r)}^2 \| H \|_{2, \B(r)}^2 dt \r)^{\frac12} \leq\\
\leq \| H \|_{2,\infty,\Q(r)} \l( \intl_{-r^2}^0 \| v \|_{3,\B(r)}^2 dt \r)^{\frac12} \leq r^{\frac32} C(r) A_*(r).
\endgathered
\label{lh31}
\end{equation}
Combining inequalities \eqref{F-4}-\eqref{lh31} we will get the statement of lemma.

\qed

Using interpolation inequality \eqref{C_3} for $C(r)$ in the right hand side \eqref{decHC}, we will obtain inequality \eqref{decHC} in another form
\begin{corollary}
Assume that conditions (\ref{Class}) --- (\ref{Divergent-free}) hold. Then for any
$0< r \leq 1$ and $0 < \tht \leq 1$ the following estimate holds
\begin{equation}
F_2(\tht r) \leq c \tht^{\al} F_2(r) + c \tht^{-\frac32} E^{\frac12}(r) A^{\frac12}(r) A_*(r).
\label{decHEA}
\end{equation}
\end{corollary}


\section{Boundedness of energy functionals}
\label{Bound_section}

In this section we derive  estimates of energy functionals which allow us to obtain uniform boundedness  (with respect to the radius)
of all functionals \eqref{EF} if boundedness of some of them is known.

Observe that one can prove a group of estimates that are the consequences of H\"older inequality, embedding theorem and interpolation inequality.
\begin{equation}\label{C_3}
C(r) \ \le  \ ~A^{\frac 12}(r)E^{\frac 12}(r),
\qquad
F(r)\ \le \ A_*^{\frac 12}(r)[ E_*^{\frac 12}(r) +  F_2^{\frac 12}(r)]
\end{equation}
\begin{equation}
D(r)\ \le \  c D_1(r), \qquad D_1(r)\ \le \ c D_s(r), \qquad \forall s>1.
\label{D}
\end{equation}

First of all we will prove the decay estimate for the pressure
\begin{lemma}
If $v,p,H$ are the suitable weak solution near the boundary to the MHD equations in $\Q$. Then for any
$0< r \leq 1$ и $0 < \tht \leq 1$ the following estimate holds
\begin{equation}
\gathered
D_{\frac{12}{11}}(\tht r) \leq c \tht^{\al}\l( D_{\frac{12}{11}}(r) + E(r) \r) +\\
+c(\tht) \l( E(r) A^{\frac12}(r) C^{\frac12}(r) + E_*(r) A_*^{\frac12}(r) F^{\frac12}(r)\r).
\endgathered
\label{D_*}
\end{equation}
\end{lemma}

\bigskip
\noindent
{\bf Proof.} To obtain (\ref{D_*}) we apply the method developed in \cite{Seregin_JMFM}, \cite{Seregin_ZNS271},
see also \cite{SSS}. Denote $\Pi_r=\mathbb R^3_+\times (-r^2,0)$. We fix $r\in (0,1]$ and $\theta\in (0,\frac 12)$
 and define a function $g:\Pi_{r}^+\to \mathbb R^3$ by the formula
$$
g \ = \  \left\{\begin{array}{cl}
 \rot H \times H-(v\cdot \nabla )v , & \mbox{in } \  Q^+(r), \\
0, & \mbox{in } \ \Pi_r^+\setminus Q^+(r)
\end{array}
\right.
$$
Then we decompose $v$ and $p$ as
$$
v\ = \ \hat v + \tilde v,\qquad p \ = \ \hat p +\tilde p,
$$
where $(\hat v, \hat p)$ is a solution of the Stokes initial boundary value problem in a half-space
$$
\gathered
\left\{\begin{array}c  \d_t \hat v - \Delta \hat v+\nabla \hat p \ = \ g, \\ \div \hat v =0 \end{array}\right. \qquad \mbox{in}\quad \Pi_r^+, \\
\hat v|_{t=0}=0, \qquad \hat v|_{x_3=0}=0,
\endgathered
$$
and $(\tilde v, \tilde p)$ is a solution of the homogeneous Stokes system in $Q^+(r)$:
$$
\gathered
\left\{\begin{array}c  \d_t \tilde v - \Delta \tilde v +\nabla \tilde p\ = \ 0, \\ \div \tilde v =0 \end{array}\right. \qquad \mbox{in}\quad Q^+(r),  \\
\tilde v|_{x_3=0}=0.
\endgathered
$$

For $\nabla \hat p$ and $\nabla \tilde p$ the following estimates hold (see \cite{Seregin_ZNS271}, see also \cite{Solonnikov_ZNS288}):
\begin{equation}
\gathered
\| \nabla \hat p\|_{L_{\frac{12}{11},\frac 32}(Q^+(r))} \ + \ \frac 1r \| \nabla \hat v \|_{L_{\frac{12}{11},\frac 32}(Q^+(r))}
  \ \le \\ \le \  c~\Big( ~\| H\times \rot H\|_{L_{\frac{12}{11}, \frac 32}(Q^+(r))} \ + \ \| (v\cdot \nabla)v \|_{L_{\frac{12}{11}, \frac 32}(Q^+(r))}~\Big),
\endgathered
\label{lh32}
\end{equation}
\begin{equation}
\| \nabla \tilde p\|_{L_{\frac{12}{11},\frac 32}(Q^+(\theta r))} \
\le \ c~\theta^{\al} ~\Big( ~\frac 1r \| \nabla \tilde v \|_{L_{\frac{12}{11},\frac 32}(Q^+(r))}
\ + \  \|  \nabla \tilde p \|_{L_{\frac{12}{11},\frac 32}(Q^+(r))}~\Big).
\label{lh33}
\end{equation}
To estimate the right hand side of \eqref{lh32} we will use H\"older and interpolation inequalities
\begin{equation}
\gathered
 \| (v \cdot \gr) v \|_{\frac{12}{11},\frac32} =\\
 = \l( \intl_{-r^2}^0 \| (v \cdot \gr) v \|_{\frac{12}{11}}^{\frac32} dt \r)^{\frac23}
 \leq \l( \intl_{-r^2}^0 \| \gr v\|_{2}^{\frac32} \| v \|_{\frac{12}{5}}^{\frac32} dt \r)^{\frac23} \leq \\
 \leq \| \gr v \|_2 \l( \intl_{-r^2}^0 \| v \|_{\frac{12}5}^6 dt \r)^{\frac16}
 \leq \| \gr v \|_2 \l( \intl_{-r^2}^0 \| v \|_{2}^3 \| v \|_3^3 dt \r)^{\frac16} \leq\\
 \leq \| \gr v\|_2 \| v \|_{2,\infty}^{\frac12} \| v \|_3^{\frac12}
\endgathered
\label{lh34}
\end{equation}
Term $\| H\times \rot H\|_{L_{\frac{12}{11}, \frac 32}(Q^+(r))}$ can be estimated similarly.

The right hand side of \eqref{lh33} can be estimated as follows
\begin{equation}
\gathered
\frac 1r \| \nabla \tilde v \|_{L_{\frac{12}{11},\frac 32}(Q^+(r))}
\ + \  \|  \nabla \tilde p \|_{L_{\frac{12}{11},\frac 32}(Q^+(r))} \leq\\
\leq c \l( \| \nabla v \|_{2,\Q(r)} + \| \nabla p \|_{L_{\frac{12}{11},\frac 32}(Q^+(r))} + \r.\\
+ \l. \frac 1r \| \nabla \hat v \|_{L_{\frac{12}{11},\frac 32}(Q^+(r))}
\ + \  \|  \nabla \hat p \|_{L_{\frac{12}{11},\frac 32}(Q^+(r))}\r).
\endgathered
\label{lh35}
\end{equation}
Combining inequalities \eqref{lh32}-\eqref{lh35} we will get the statement of lemma.

\qed


\begin{theorem}
Let $v,p,H$ are the suitable weak solution near the boundary to the MHD equations in $\Q$ and
$$C(R) + F_2(R) \leq M, \quad 0 < R \leq 1.$$
If we consider the following functional
$$\L(r) = A^2(r) + E^2(r) + A_*^2(r) + E_*^2(r) + D_{\frac{12}{11}}^{\frac{25}{24}}(r),$$
then the following estimate will hold
$$\L(r) \leq C(M) (r^{\al} \L(1) + 1).$$
\end{theorem}

\bigskip
\noindent
{\bf Proof.} From local energy inequality we obtain
$$\L( \tht r) \leq c \l( C_2^2(2 \tht r) + F_2^2(2 \tht r) + C^3(2 \tht r) + C(2 \tht r) D(2 \tht r) + \r.$$
$$\l. C^2(2 \tht r)F_3 (2 \tht r) + C(2 \tht r)F_3^2 (2 \tht r)+ D_{\frac{12}{11}}^{\frac{25}{24}}(\tht r) \r).$$
Now we will estimate every term in the right hand side. Our goal is to prove the following estimate
\begin{equation}
\L(\tht r) \leq \frac12 \L(r) + C(M).
\label{lb2}
\end{equation}
Then we can use a standard iteration procedure (see \cite{Ser1}) and obtain the statement of the theorem.

Estimates for the first three terms are obvious. To estimate the 4th term we use Young inequality and inequality \eqref{D}
$$ C(2 \tht r) D(2 \tht r) \leq c \l( D_{\frac{12}{11}}^{\frac{25}{24}}(2 \tht r) + M^{25}\r). $$
Now we are going to prove estimate for $D_{\frac{12}{11}}(2 \tht r)$. To do this, we will use inequalities \eqref{C_3} and \eqref{D_*}
\begin{equation}
\gathered
D_{\frac{12}{11}} (2 \tht r) \leq \\ \leq c \tht^{\al} \l( D_{\frac{12}{11}}(r) + E(r)\r) +
c(\tht) \l( E(r) A^{\frac12}(r) C^{\frac12}(r) + E_*(r) A_*^{\frac12}(r) F_3^{\frac12}(r) \r) \leq \\
\leq c \tht^{\al} \l( D_{\frac{12}{11}}(r) + E(r)\r) + c(\tht) \l( \L^{\frac34}(r) M^{\frac12} + \L^{\frac34}(r) F_3^{\frac12}(r)\r).
\endgathered
\label{lb1}
\end{equation}
We use interpolation inequality to estimate $F_3(r)$
\begin{equation}
\gathered
F_3(r) \leq F_{\frac{10}{3}}^{\frac56}(r) F_2^{\frac16}(r)
\leq c \l(A_*^{\frac25}(r) \l( E_*^{\frac35}(r) + F_2^{\frac35}(r)\r)\r)^{\frac56} F_2^{\frac16}(r) \leq \\
\leq c \l( \L^{\frac5{12}}(r) M^{\frac16} + \L^{\frac16}M^{\frac23}\r).
\endgathered
\label{lb3}
\end{equation}
Now we substitute this to \eqref{lb1}
$$
\gathered
D_{\frac{12}{11}}( 2 \tht r) \leq c \tht^{\al} \l( D_{\frac{12}{11}}(r) + E(r)\r) +\\
+ c(\tht) \l( \L^{\frac34}(r) M^{\frac12} + \L^{\frac{23}{24}}(r) M^{\frac{1}{12}}
 + \L^{\frac56}(r) M^{\frac13}\r).
\endgathered
$$
As the result we obtain
$$ D_{\frac{12}{11}}^{\frac{25}{24}} (2 \tht r) \leq c \tht^{\al} \L(r) + c(\tht) \l( \L^{\frac{575}{576}}(r) M^{k_1} + \L^{\frac{125}{144}}(r) M^{k_2} \r).$$
Since the right hand side of the last inequality contain $\L(r)$ in the degree smaller then 1, choosing $\tht$ sufficiently small and using Young inequality
we obtain an estimate \eqref{lb2}.

To estimate the last two terms we use \eqref{lb3}
$$ C(r)F_3^2(r) \leq c \l( \L^{\frac56}(r) M^{\frac43} + \L^{\frac13} M^{\frac73} \r)$$
and Young inequality. The second term can be estimated in the same maner.

As the result we obtain \eqref{lb2}. Next by standard iteration procedure we finish the prove of the theorem.

\qed

\section{Proof of main results}
\label{reg_section}

As a first step we obtain theorem \ref{Main_ep-reg} without smallness condition on a pressure.
\begin{lemma}
\label{MainSmallLemma}
For arbitrary $M>0$ there is $\ep_1(M)>0$, such that if $v,p,H$ are the suitable weak solution near the boundary to the MHD equations in $\Q$,
\begin{equation}
A(R) + E(R) + A_*(R) + E_*(R) + F_3(R) + D_{\frac{12}{11}}(R) < M, \quad \forall \, 0 < R \leq 1
\label{ls1}
\end{equation}
and
\begin{equation}
C(1) + F_3(1) < \ep_1,
\label{ls2}
\end{equation}
then the functions $v$ and $H$ are H\" older continuous on $\bar Q^+(r_*)$ for some $0 < r_* < 1$.
\end{lemma}

\bigskip
\noindent
{\bf Proof.} Assume that the statement of the lemma is false. Then there are sequences of $v_n,p_n,H_n$ of suitable weak solutions in $\Q$,
such that
\begin{equation}
C(v_n,1) + F_3(H_n,1) = \ep_n \to 0, \text{ as } n \to \infty
\label{ls3}
\end{equation}
and $0$ is a singular point. Then by theorem \ref{Main_ep-reg}
\begin{equation}
C(v_n,r) + D(p_n,r) + F_3(H_n,r) > \ep_*
\label{ls4}
\end{equation}
for all $0 < r < 1$.

On the other hand from \eqref{D_*}, \eqref{ls1}, \eqref{ls3} and the embedding theorem we have
\begin{equation}
D(p_n,r) \leq c D_{\frac{12}{11}}(p_n, r) \leq c r^{\al} M + c(r) M^{\frac32} \ep_n^{\frac12}.
\label{ls5}
\end{equation}
So we fix $0 < r \leq 1$ and pass to the limit by $n$ in \eqref{ls4} and \eqref{ls5}
$$
\gathered
 \ep_* \leq \limsup_{n \to \infty} \l( C(v_n,r) + D(p_n,r) + F_3(H_n,r) \r) = \\
= \limsup_{n \to \infty} D(p_n,r) \leq c r^{\al} M.
\endgathered
$$
As the result we obtain, that the inequality
$$ \ep_* \leq c r^{\al} M$$
must be true for arbitrary $0 < r < 1$. So we have a contradiction.

\qed


\begin{theorem}
For arbitrary $M>0$ there is $\ep_2(M)>0$, such that if $v,p,H$ are the suitable weak solution near the boundary to the MHD equations in $\Q$,
satisfying to \eqref{ls1} and one of the following conditions holds
\begin{equation}
E(1) < \ep_2,
\label{Econd}
\end{equation}
\begin{equation}
A(1) < \ep_2,
\label{Acond}
\end{equation}
\begin{equation}
C(1) < \ep_2,
\label{Ccond}
\end{equation}
then the functions $v$ and $H$ are H\" older continuous on $\bar Q^+(r_*)$ for some $0 < r_* < 1$.
\label{Ecrit}
\end{theorem}

\bigskip
\noindent
{\bf Proof.} The proof of this theorem is similar to the proof of previous lemma. We begin from the case \eqref{Econd}.
Let $v_n,p_n,H_n$ are the sequences of suitable weak solutions to the MHD system, such that \eqref{ls1} holds,
$$E(v_n,1) = \ep_n \to 0$$
as $n \to \infty$, and $z_0=0$ is a singular point. Then from lemma \ref{MainSmallLemma} we have
\begin{equation}
C(v_n,r) + F_3(v_n,r) > \ep_1
\label{ls9}
\end{equation}
for arbitrary $0 < r < 1$.

On the other hand
\begin{equation}
C(v_n,r) \leq \frac{c}{r^{\frac23}} C(v_n,1) \leq \frac{c}{r^{\frac23}} A^{\frac12}(v_n,1) E^{\frac12}(v_n,1) \to 0
\label{ls7}
\end{equation}
as $n \to \infty$ and for any fixed $0< r \leq 1$. From \eqref{decHEA} we obtain
\begin{equation}
\limsup F_2(H_n,r) \leq c r^{\al} M.
\label{ls8}
\end{equation}
Next we use interpolation inequality
\begin{equation}
F_3(H_n,r) \leq F_2^{\frac16}(H_n,r) F_{\frac{10}{3}}^{\frac56}(H_n,r).
\label{ls41}
\end{equation}
To estimate the second factor in the right hand side of \eqref{ls41} we use \eqref{C_3}. So from \eqref{ls9}-\eqref{ls41} we obtain
$$ \ep_1 \leq \limsup_{n \to \infty} \l(C(v_n,r) + F_3(v_n,r) \r) \leq c r^{\al_1} M^k \quad \forall \, 0 < r \leq \frac12,$$
and, if we choose $r$ sufficiently small, we will have a contradiction.

Observe, that $E(r)$ and $A(r)$ take part in \eqref{decHEA} symmetrically, so the proof of this theorem in the case \eqref{Acond} is similar to the previous one.
In the case of \eqref{Ccond} for obtaining \eqref{ls8} is sufficient to use \eqref{decHC}.

\qed

\end{document}